\documentclass[12pt]{article}
\usepackage{cite}
\usepackage{amsmath,amssymb,amsfonts}
\usepackage{algorithm}
\usepackage{graphicx}

\newtheorem{proposition}{Proposition}
\newtheorem{theorem}{Theorem}

\newtheorem{remark}{Remark}
\newtheorem{lemma}{Lemma}

\newtheorem{definition}{Definition}

\newcommand{\defeq}{\overset{\mathrm{def}}{=}}
\newcommand{\tr}{\top}
\newcommand{\qed}{\hfill $\square$}

\begin{document}

\title{Parallel MPC for Linear Systems with State and Input Constraints}

\author{Jiahe Shi\thanks{Jiahe Shi and Boris Houska are with ShanghaiTech University, China. (email: {\tt{shijh, borish}@shanghaitech.edu.cn})} , Yuning Jiang\thanks{Yuning Jiang is with Automatic Control Laboratory, EPFL, Switzerland. (email: {\tt yuning.jiang@epfl.ch})} , Juraj Oravec\thanks{JO gratefully acknowledges financial support of the Scientific Grant Agency of the Slovak Republic under the grant VEGA 1/0297/22.} \thanks{Juraj Oravec is with Faculty of Chemical and Food Technology, Slovak University of Technology in Bratislava, Slovakia. (email: \tt juraj.oravec@stuba.sk)
} , and Boris Houska
}

\date{}

\maketitle

\thispagestyle{empty}

\begin{abstract}
This paper proposes a parallelizable algorithm for linear-quadratic model predictive control (MPC) problems with state and input constraints. The algorithm itself is based on a  parallel MPC scheme that has originally been designed for systems with input constraints. In this context, one contribution of this paper is the construction of time-varying yet separable constraint margins ensuring recursive feasibility and asymptotic stability of sub-optimal parallel MPC in a general setting, which also includes state constraints.  Moreover, it is shown how to tradeoff online run-time guarantees versus the conservatism that is introduced by the tightened state constraints. The corresponding performance of the proposed method as well as the cost of the recursive feasibility guarantees is analyzed in the context of controlling a large-scale mechatronic system. This is illustrated by numerical experiments for a large-scale control system with more than $100$ states and $60$ control inputs leading to run-times in the millisecond range.\\[0.1cm]

\end{abstract}

%\begin{keywords}
%
%\end{keywords}

\section{Introduction}
\label{sec:introduction}

Model predictive control (MPC)~\cite{Rawlings2009} is a modern optimization based control technique. In many industrial applications of MPC~\cite{Qin2003} linear models with quadratic costs are used, such that the online optimization problems can be formulated as quadratic programming (QP) problems. As state-of-the-art centralized solvers, including active-set solvers~\cite{Ferreau2014} and real-time interior point solvers~\cite{FRISON2020}, can solve moderately sized QPs within the milli- to microsecond range, the run-time of these solvers is hardly ever a problem for medium-scaled MPC problems. For larger problems, however, one needs to use first order methods, which can exploit the sparsity and structure of the online QPs. Examples for such first order methods include dual decomposition~\cite{Everett1963}, ADMM~\cite{Boyd2011}, and ALADIN~\cite{Houska2016} with constant Hessian approximations, which have all been further adapted for solving distributed MPC problems~\cite{Giselesson2013,Conte2012,Necoara2008,Donoghue2013}. In practice, first order methods often need thousands of iterations until a sufficient numerical accuracy is achieved. This can be a problem for some large-scale applications. Moreover, if one is interested in providing stability guarantees, it is not always clear how to trade off the run-time versus the numerical accuracy of these solvers.

\bigskip
\noindent
The present paper proposes a non-trivial extension of the parallel MPC scheme from~\cite{Jiang2021}, which combines ideas from both first order methods and Explicit MPC methods~\cite{BemEtal:aut:02}. In this context, it is important to understand first that Explicit MPC can usually only be used for systems with a limited number of constraints such that one can solve associated multi-parametric QPs offline, by pre-computing piecewise affine (PWA) solution maps~\cite{BemEtal:aut:02}. These PWA maps are then evaluated online for real-time control~\cite{Ingole2015}. However, the number of critical regions of PWA solution maps of MPC problems grows, in the worst case, exponentially with the number of constraints~\cite{borrelli2003geometric}. This renders traditional Explicit MPC essentially inapplicable to large-scale MPC. Nevertheless, if one uses a parallel MPC scheme, one can solve the smaller-scale distributed QPs by using methods from Explicit MPC~\cite{Jiang2021}.

\bigskip
\noindent
Apart from the above literature on QP solvers for MPC, much literature can be found on the stability and recursive feasibility of distributed MPC controllers in the presence of state constraints. For instance, in~\cite{FERRAMOSCA2013} a terminal set for cooperative control is designed in order to track changing set-points. In~\cite{Conte2016}, a separable terminal cost combined with time-varying local terminal sets is introduced. And, in~\cite{Darivianakis::TAC19}, an adaptive terminal region computation scheme is used to reduce conservatism. Comprehensive reviews of the application of invariant (or contractive) sets in the context of stability guarantees for distributed MPC can be found in~\cite{HERNANDEZ201711829,Cannon::TAC::2017::7488973,Matthias17}.

\subsection*{Outline}
Section~\ref{sec::parallel_explicit_mpc} reviews existing parallel MPC methods for systems with input constraints~\cite{Jiang2021}. The main theoretical contribution of the current article is presented in Section~\ref{sec::contractive_sets_and_feasibility}, which discusses a general strategy for constructing separable time-varying robustness margins for state-constrained linear systems by using contractive ellipsoidal sets. Section~\ref{sec::fmpc} explains how these time-varying constraints can be featured within the context of real-time parallel MPC while maintaining both recursive feasibility as well as asymptotic closed-loop stability of the associated sub-optimal MPC controller. Section~\ref{sec::case_study} presents a numerical case study for a large-scale MPC problem with more than $100$ system states.

\subsection*{Notation} 
We use the notation $\mathrm{diag}(v)$ to denote the diagonal matrix in $\mathbb R^{n \times n}$, whose diagonal elements are the coefficients of a vector $v \in \mathbb R^{n}$. The Minkowski sum and Pontryagin difference of two given sets $X,Y \subseteq \mathbb R^n$ are denoted by
\begin{equation}
\begin{aligned}
X \oplus Y &\; \defeq \; \left\{ x+y \ \middle| \ x \in X, \ y \in Y \right\}  \\
\text{and} \qquad X \ominus Y &\; \defeq \; \left\{ z  \ \middle| \ \{ z \} \oplus Y \subseteq X \right\} .
\end{aligned}
\end{equation}

\section{Parallel Explicit MPC}
\label{sec::parallel_explicit_mpc}
This paper concerns the linear-quadratic MPC problem
\begin{eqnarray}
J_N(\hat x) \ \defeq &\underset{x,u}{\min}&\;\;x_N^\tr P x_N+\sum_{k=0}^{N-1} \left( x_k^\tr Q x_k + u_k^\tr R u_k \right) \notag \\
&\mathrm{s.t.}& \label{eq::MPC}
\;\left\{
\begin{aligned}
&\forall k \in \{ 0, \ldots, N-1 \}, \\
&x_{k+1} = A x_k + B u_k \quad \mid \, \lambda_k, \\
&x_0 = \hat x, \ (x_k, u_k) \in \mathbb{X} \times \mathbb{U} \; ,
\end{aligned}
\right.
\end{eqnarray}
where $A$ and $B$ denote system matrices of appropriate dimension, $Q$, $R$, $P$ positive definite tuning parameters, $x_k$ and $u_k$ the states and controls, and $\lambda_k$ the co-states of~\eqref{eq::MPC}. Moreover, $\hat x \in \mathbb{R}^{n_{\mathrm{x}}}$ denotes the initial measurement. Polyhedral state and control constraints are given by
\begin{subequations}
\begin{align}
\label{eq::defX}
\mathbb X & = \{ x \in \mathbb R^{n_x} \mid Cx \leq c  \} \\
\label{eq::defU}
\text{and} \qquad \mathbb U &= \{ u \in \mathbb R^{n_u} \mid Du \leq d  \},
\end{align}
\end{subequations}
where $C$ and $D$ are given matrices and $c > 0$, $d > 0$ are associated bounds, such that the point $(0,0) \in \mathbb X \times \mathbb U$ is strictly feasible. The latter condition ensures that $P$ can be found by solving an algebraic Riccati equation~\cite{Rawlings2009}, such that we have $J_N(\hat x) = J_\infty(\hat x)$ for any sufficiently large prediction horizon $N$ and any initial state $\hat x$.
\begin{remark}
The assumption that $N$ is sufficiently large, such that $J_N = J_\infty$, is introduced for simplicity of presentation. Alternatively, one can leave the terminal cost away and rely on turnpike horizon bounds~\cite{Damm2014}.
\end{remark}

\subsection{Parallel MPC with input constraints}
\label{sec::pmpc_woth_input_constraints}
This section reviews an existing parallel MPC scheme from~\cite{Jiang2021} for systems with input constraints but $\mathbb{X} = \mathbb{R}^{n_x}$, where~\eqref{eq::MPC} is solved by starting with initial guesses
\[
z = [z_0^\top,z_1^\top, \ldots, z_N^\top ]^\top \quad \text{and} \quad v = [v_0^\top,v_1^\top, \ldots, v_{N-1}^\top ]^\top
\]
for the state and control trajectories, as well as an initial guess for the co-state $\lambda =[\lambda_0^\top, \lambda_1^\top, \ldots, \lambda_{N-1}^\top]^\top$.
The method proceeds by solving the optimization problems\footnote{Notice that~\eqref{eq::initial},~\eqref{eq::horizon}, and~\eqref{eq::terminal} are parametric QPs for which it is sometimes possible to precompute explicit solution maps for online evaluation~\cite{BorrelliPHD,MPT3}. In particular, if $Q,R,P$ in~\eqref{eq::MPC} and $D$ in~\eqref{eq::defU} are block-diagonal, each of the decoupled QPs is itself separable and thus can be further parallelized.}
\begin{align}
\label{eq::initial}
\underset{u_0 \in \mathbb U}{\min} \quad &  \Vert u_0 \Vert_R^2 - \lambda_0^\tr B u_0 + \| u_0 - v_0 \|_R^2 \, ,
\end{align}
\begin{align}
\underset{x_k,u_k \in \mathbb U}{\min} \;\; &\Vert x_k \Vert_Q^2 + \Vert u_k \Vert_R^2- \lambda_k^\tr B u_k + (\lambda_{k-1} - A^\tr \lambda_k)^\tr x_k  \notag \\
\label{eq::horizon}
& + \| x_k - z_k \|_Q^2 + \| u_k - v_k \|_R^2 \; ,
\end{align}
\begin{align}
\label{eq::terminal}
\text{and} \qquad \underset{x_{N}}{\min} \;  \; &\Vert x_{N} \Vert_P^2 + \lambda_{N-1}^\tr x_{N} + \| x_N - z_N \|_P^2 
\end{align}
in parallel with $k \in \{ 1, \ldots, N-1 \}$ in~\eqref{eq::horizon}. If $(z,v,\lambda)$ is an optimal solution of~\eqref{eq::MPC}, the solutions $x_k$ and $u_k$ of~\eqref{eq::initial},~\eqref{eq::horizon}, and~\eqref{eq::terminal} are an optimal solution of~\eqref{eq::MPC}. However, in general, our initial guesses are not optimal. Therefore, the current iterates for $x$ and $u$ might not even correspond to a feasible trajectory. Consequently, we solve a consensus problem\footnote{The consensus QP~\eqref{eq::consensusQP} can be viewed as a parametric LQR problem, for which a linear explicit solution map can be pre-computed~\cite{Bertsekas2012}.} to update the variables $(z,v,\lambda)$
\begin{align}
\min_{z^+,v^+} \; &
\sum_{k=0}^{N-1} \left\|
\begin{bmatrix}
z_k^+ + z_k - 2 x_k \\[0.12cm]
v_k^+ + v_k - 2 u_k
\end{bmatrix}
\right\|_{\Sigma}^2 +  \Vert z_N^+ + z_N - 2 x_N \Vert_P^2 \notag \\[0.12cm]
\label{eq::consensusQP}
\text{s.t.}\; \;  & \; \left\{
\begin{aligned}
&\forall k \in \{ 0, \ldots, N-1 \}, \\[0.12cm]
&z_{k+1}^+ = A z_k^+ + B v_k^+ \mid \; \delta_k, \\[0.12cm]
&z_0^+ = \hat x
\end{aligned}
\right.
\end{align}
with $\Sigma = \text{diag}(Q,R)$.
The decoupled QPs~\eqref{eq::initial},~\eqref{eq::horizon}, and~\eqref{eq::terminal} and the consensus QP~\eqref{eq::consensusQP} need to be solved repeatedly, in an alternating way, as detailed in Algorithm~\ref{alg::mpc}.

\begin{algorithm}[htbp!]
\small
\caption{Parallel MPC with input constraints.}
\textbf{Initialization:}
Guesses for the states, inputs, and co-states $(z,v,\lambda)$.

\textbf{Online:}
\begin{enumerate}

\item For $m = 1: \overline{m}$:

\begin{enumerate}

\item Solve the decoupled QPs~\eqref{eq::initial},~\eqref{eq::horizon}, and~\eqref{eq::terminal} in parallel to compute the optimal solution
$$x = [x_1,\ldots,x_N] \quad \text{and} \quad u = [u_0,\ldots,u_{N-1}] \; .$$

\item Compute an optimal primal-dual solution $(z^+,v^+)$ of the consensus QP~\eqref{eq::consensusQP}, and update
$$z \leftarrow z^+, \ v \leftarrow v^+, \ \lambda \leftarrow \lambda + \delta\; .$$
\end{enumerate}

\item Send $u_0$ to the real plant.

\item Shift all variables and go to Step 1,
\begin{align}
z &\leftarrow [z_1, \ldots, z_N, 0]\;, \quad v \leftarrow [v_1, \ldots, v_{N-1}, 0]\;, \notag \\
\lambda &\leftarrow [\lambda_1, \ldots, \lambda_{N-1}, 0] \; .\notag 
\end{align}
\end{enumerate}
\label{alg::mpc}
\end{algorithm}

\bigskip
\noindent
Because Step~1) of this algorithm implements a finite number of iterations $\overline{m}$ at each MPC step, Step~3) sends only an approximation, $u_0 \approx u_0^\star$, of the optimal control input $u_0^\star$ to the real plant. Therefore, there arises the question in which sense Algorithm~\ref{alg::mpc} yields a stable, let alone feasible controller.

\subsection{Convergence rate estimates}
\label{sec::convergence}
In order to review the convergence properties of Algorithm~1, we introduce the auxiliary function
\begin{align}
\Phi(z,v,&\lambda) \defeq  \sum_{k=0}^{N-1} \left( z_k^\tr Q z_k + v_k^\tr R v_k \right) + z_N^\tr P z_N \notag \\
& + \frac{1}{4} \left( \Vert B^\tr \lambda_0 \Vert_{R^{-1}}^2 +  \sum_{k=1}^{N-1} \Vert B^\tr \lambda_k \Vert_{R^{-1}}^2 \right) \notag \\
& + \frac{1}{4} \left( \sum_{k=1}^{N-1} \Vert \lambda_{k-1}-A^\tr \lambda^k \Vert_{Q^{-1}}^2 +  \Vert \lambda_{N-1} \Vert_{P^{-1}}^2 \right),
\notag
\end{align}
which corresponds to the sum of the objective function of~\eqref{eq::MPC} and its weighted conjugate \mbox{function---see}~\cite[Sect.~II.B]{Jiang2021} for details about the properties of this function. Notice that $\Phi$ is a positive definite quadratic form that can be used to measure the distance $\Delta$ from the current iterate to the primal-dual solution $(x^\star,u^\star,\lambda^\star)$ of~\eqref{eq::MPC}, given by
\begin{equation}
\label{def::distanceToOptimalSol}
\Delta(z,v,\lambda) \,\defeq \, \Phi( z - x^\star, v - u^\star, \lambda - \lambda^\star ) \; .
\end{equation}
A proof of the following lemma follows by combining the results from~\cite[Thm.~1 \& Thm.~2]{Jiang2021}.

\begin{lemma}
\label{lem::distance_to_optimal_finite_iter}
Let $u_0$ denote the approximately optimal input that is sent to the real process in Step~2 of Algorithm~\ref{alg::mpc} after running a finite number of iterations $\overline m$. Moreover, let
\[
x_0^+ \, \defeq \, A \hat x + B u_0 \; ,
\]
denote the approximately optimal closed-loop state at the next time instance and $x_1^\star$ the optimal state that would be reached if the optimal input $u_0^\star$ would be send to the real plant. Then there exists a constant $\kappa < 1$ such that
\begin{eqnarray}
\label{eq::MAIN_INEQUALITY}
\Vert x_0^+ - x_1^\star \Vert \leq \sigma (1+\kappa) \kappa^{\overline{m}+1} \Delta(z,v,\lambda) ,
\end{eqnarray}
for any constant $\sigma > 0$ that satisfies $B^T Q B \preceq \sigma R$, where $z,v,\lambda$ denotes the initialization of Algorithm~\ref{alg::mpc}.
\end{lemma}

\section{Contractive Sets and Feasibility}
\label{sec::contractive_sets_and_feasibility}

Throughout the following derivations, we assume that the pair $(A,B)$ is asymptotically stabilizable recalling that this is the case if and only if one can find a linear feedback gain $K$ for which the spectral radius of the matrix $A+BK$ is strictly smaller than $1$. If this assumption holds, there exists a contractivity constant $\beta \in \mathbb{R}$ satisfying
\begin{align}
\label{eq::beta}
\rho(A+BK) < \beta < 1 ,
\end{align}
where $\rho(\cdot)$ denotes the spectral radius. The following sections present various technical developments based on this assumption on $(A,B)$ that will later be needed in Section~\ref{sec::fmpc} to construct a recursively feasible parallel MPC controller.

\subsection{Ellipsoidal contractive sets}
\label{sec::ellipoids}
This section reviews the standard definition of $\beta$-contractive sets, which reads as follows.
\begin{definition}
A set $\mathbb Z \subseteq \mathbb X$ is called $\beta$-contractive if
\[
\forall x \in \mathbb Z, \, \exists u \in \mathbb U, \quad A x + Bu \in \beta \mathbb Z \; .
\]
\end{definition}
Notice that $\beta$-contractive sets exist whenever~\eqref{eq::beta} holds. In particular, a $\beta$-contractive ellipsoid of the form
$$\mathbb Z = \mathcal E(Z) \; \defeq \; \left\{ Z^\frac{1}{2}v \, \middle| \, \Vert v \Vert_2^2 \leq 1 \,  \right\} , $$
with shape matrix $Z$ can be found by solving the convex semi-definite programming (SDP) problem
\begin{equation}
\label{eq::maxZ}
\begin{aligned}
 \min_{Z} \quad &\;\;\mathrm{tr}(Z) \\ \text{s.t.}\quad &  \left\{
\begin{array}{l}
(A+BK) Z (A+BK)^\tr \; \preceq \; \beta^2 Z ,  \\[0.16cm]
r^2 I \; \preceq \; Z , \\[0.16cm]
Z \; = \; Z^\tr \succ 0 , \\[0.16cm]
C Z C^\tr \; \leq \; (1+\alpha)^{-2} \cdot \mathrm{diag}^2(c) \\[0.16cm]
D K Z K^\tr D^\tr \; \leq \; (1+\alpha)^{-2} \cdot \mathrm{diag}^2(d) 
\end{array}
\right. 
\end{aligned}
\end{equation}
for a given constant $\alpha \geq \beta^N$ and a given inner radius $r > 0$.
\begin{proposition}
\label{prop::ellipsoid}
Let $\mathbb X$ and $\mathbb U$ be given as in~\eqref{eq::defX} and~\eqref{eq::defU} with $c > 0$ and $d > 0$, and let $\beta$ satisfy~\eqref{eq::beta}. If the radius $r > 0$ is sufficiently small, then~\eqref{eq::maxZ} admits a minimizer $Z$ and $\mathcal E(Z)$ is a $\beta$-invariant set with inner radius~$r$.
\end{proposition}
\textit{Proof.}
Lyapunov inequalities of the form $S Z S^\tr \preceq Z$ admit a symmetric and positive definite solution if the eigenvalues of $S$ are all in the open unit disk. Because $\beta$ satisfies~\eqref{eq::beta},
\[
S = \beta^{-1} (A+BK)
\]
has this property and the first semi-definite inequality in~\eqref{eq::maxZ}, together with the condition $Z = Z^\tr$ has a positive definite solution. These Lyapunov conditions are equivalent to enforcing the ellipsoid $\mathcal E(Z)$ to be contractive~\cite{Blanchini2008}, while the feasibility constraints
\[
(1+\alpha) \mathcal E(Z) \subseteq \mathbb X \quad \text{and} \quad (1+\alpha) K \mathcal E(Z) \subseteq \mathbb U
\]
are equivalent to enforcing the last two inequalities in~\eqref{eq::maxZ}.
Since the Lyapunov contraction constraint is homogeneous in $Z$ while $c > 0$ and $d > 0$, we can find a small inner radius $r > 0$ of $\mathcal E(Z)$ for which all inequalities in~\eqref{eq::maxZ} are strictly feasible. As the objective of~\eqref{eq::maxZ} is bounded from below by $\mathrm{tr}(Z) \geq n_x r^2$, the statement of the proposition follows.\qed
\begin{remark}
The trace of the matrix $Z$ in the objective of~\eqref{eq::maxZ} could be replaced by other measures that attempt minimizing the size of the ellipsoid $\mathcal E(Z)$. For instance, one could also minimize the determinant or the maximum eigenvalue of $Z$.
\end{remark}

\subsection{Separable safety margins and terminal regions}
\label{sec::margins}
Because $\mathbb Z = \mathcal E(Z)$ is an ellipsoid, this set is not separable and, consequently, may not be used directly as constraint or terminal region, as this would be in conflict with our objective to exploit the separable structure of~\eqref{eq::MPC}. Nevertheless, an important observation is that the polyhedra
\begin{subequations}
\begin{align}
\label{eq::defW}
\mathbb W_k &\defeq \mathbb U \ominus (1-\beta^k) K \mathbb Z \\[0.16cm]
\label{eq::defY}
\text{and} \quad \mathbb Y_k &\defeq \mathbb X \ominus (1-\beta^k) \mathbb Z 
\end{align}
\end{subequations}
admit explicit separable representations whenever $\mathbb X$ and $\mathbb U$ are separable---despite the fact that $\mathbb Z$ is not separable. These polyhedral sets are represented as $\mathbb Y_k = \{ \, x \, \mid \, C x \leq \hat c_k \, \}$ and $\mathbb W_k = \{\, u\,\mid\, Du \leq \hat d_k \,\}$, where
\begin{subequations}
\label{eq::modified_ineq_bounds}
\begin{align}
\hat c_k \ &\defeq \ c - (1-\beta^k) \left[
\begin{smallmatrix} % \begin{array}{c}
\sqrt{C_1 Z C_1^\tr} \\[-0.1cm]
\vdots \\[0.1cm]
\sqrt{C_{n_c} Z C_{n_c}^\tr}
\end{smallmatrix} %\end{array}
\right]\\
\text{and}\quad\hat d_k \ &\defeq \
d - (1-\beta^k) \left[
\begin{smallmatrix} % \begin{array}{c}
\sqrt{D_1 KZK^\top D_1^\tr} \\[-0.15cm]
\vdots \\[0.1cm]
\sqrt{D_{n_d} KZK^\top D_{n_d}^\tr}
\end{smallmatrix} %\end{array}
\right] \; .
\end{align}
\end{subequations}
Thus, $\mathbb Y_k$ is trivially separable whenever $C$ is block-diagonal and the same holds for $\mathbb W_k$ whenever $D$ is block-diagonal. We recall that Proposition~\ref{prop::ellipsoid} ensures that the sets $\mathbb W_k$ and $\mathbb Y_k$ are non-empty if $Z$ is a feasible solution of~\eqref{eq::maxZ}.

\bigskip
\noindent
In the following, we introduce the terminal region
\begin{align}
\label{eq::defT}
\mathbb T \defeq  \alpha \mathbb Z \, ,
\end{align}
where $\alpha \geq \beta^N$ is the constant used in~\eqref{eq::maxZ}. As this terminal region is---in contrast to the sets $\mathbb Y_k$ and $\mathbb W_k$---\textit{not} separable, we will have to discuss later on how we actually avoid implementing it. However, for the following theoretical developments, it is convenient to temporarily introduce~\eqref{eq::defT}. Because $\alpha$ is used to scale the last two inequalities in~\eqref{eq::maxZ}, $\mathbb T$ is a feasible $\beta$-contractive set. It satisfies
\begin{equation}
\label{eq::Tproperty}
(A+BK) \mathbb T \oplus \beta^N (1-\beta) \mathbb Z \subseteq ( \alpha \beta + \beta^N(1-\beta)) \mathbb Z \subseteq\mathbb T ,
\end{equation}
where the first inclusion follows from the fact that $\mathbb Z$ is $\beta$-contractive, while the second follows by substituting the inequality $\alpha \geq \beta^N$. Moreover, since the factor $(1+\alpha)$ has been introduced in~\eqref{eq::maxZ}, $\mathbb T \subseteq \mathbb Y_N$ holds by construction.

\subsection{Admissible set}
\label{sec::admissible}

One of the core technical ideas of this paper is to introduce a set of admissible state and control pairs, given by
\begin{equation}
\label{def::admissibleset}
\mathcal A \defeq \left\{ 
(x,u)  \middle|
\begin{array}{l}
\exists y,w: \; \forall k \in \{ 0, \ldots, N-1 \}, \\[0.16cm]
(y_0, w_0) = (x, u), \\[0.16cm]
(y_k, w_k) \in \mathbb{Y}_k \times \mathbb{W}_k \\[0.16cm]
y_{k+1} = A y_k + B w_k , \; y_N \in \mathbb T
\end{array}
\right\} ,
\end{equation}
where the separable sets $\mathbb Y_k$ and $\mathbb W_k$ are defined as in~\eqref{eq::defW} and~\eqref{eq::defY}, respectively. We recall, the set $\mathbb T$, as defined in~\eqref{eq::defT}, is not separable but introduced temporarily for the sake of analyzing the set $\mathcal A$ in~\eqref{def::admissibleset}. The following lemma establishes a robust recursive feasibility result that turns out to be of high practical relevance for the construction of distributed MPC controllers.
\begin{lemma} \label{lem::recursively:admissible}
If $(x,u)$ is an element of the admissible set, $(x,u) \in \mathcal A$, then there exists for every $x^+$ satisfying 
\begin{align}
\label{eq::perturbation}
\left\| x^+ - (A x+B u) \right\|_2 \; \leq \; (1-\beta) \, r
\end{align}
a control $u^+$ such that $(x^+,u^+) \in \mathcal A$.
\end{lemma}

\textit{Proof.} Let the pair $(y,w)$ satisfy the conditions in the above definition of $\mathcal A$ for the pair $(x,u) \in \mathcal A$. This implies
\[
y_1 = Ax + Bu \in \mathbb Y_1 = \mathbb X \ominus (1-\beta) \mathbb Z
\]
and, consequently, if $x^+$ satisfies~\eqref{eq::perturbation}, we must have $x^+ \in \mathbb X$, since $r$ is an inner radius of $\mathbb Z$ by construction. Thus, we define the shifted initial value
\begin{eqnarray}
\label{eq::inductionStart}
y_0^+ \, \defeq \, x^+ \in \mathbb Y_0 \; .
\end{eqnarray}
Next, a feasible trajectory is generated via the closed-loop recursion
\begin{eqnarray}
\label{eq::y+}
y_{k+1}^+ &=& y_{k+2} + (A+BK)(y_k^+-y_{k+1}) , \qquad \\[0.1cm]
\label{eq::w+}
w_k^+ &=& u_{k+1} + K( y_k^+ - y_{k+1} ) ,
\end{eqnarray}
for $k \in \{ 0, \ldots, N-1 \}$, where we additionally define
\begin{align}
\label{eq::terminalShift}
y_{N+1} \; \defeq \; (A+BK) y_N \qquad \text{and} \qquad u_{N} \; \defeq \; K y_N ,
\end{align}
such that~\eqref{eq::y+} and~\eqref{eq::w+} are well-defined for all $k$, including the special case $k=N-1$. In order to check that $(y^+,w^+)$ satisfies the system dynamic, we briefly verify that
\begin{eqnarray}
y_{k+1}^+ &\overset{\eqref{eq::y+},\eqref{eq::terminalShift}}{=}&  A y_{k+1} + B u_{k+1} + (A+BK)(y_k^+ - y_{k+1})  \notag \\[0.16cm]
&=& A y_k^+ + B( u_{k+1} + K(y_k^+ - y_{k+1})) \notag \\[0.16cm]
&\overset{\eqref{eq::w+}}{=}& A y_k^+ + B w_k^+ \; .
\end{eqnarray}
Next, we use an induction to show that
\begin{align}
\label{eq::InductionClaim}
y_k^+ - y_{k+1} \in \beta^k (1-\beta) \mathbb Z 
\end{align}
holds for all $k \in \{ 0,1,\ldots,N \}$. For $k=0$, this follows from~\eqref{eq::inductionStart}, since $y_0^+ = x^+$ and $x^+ - y_1 \in (1-\beta)\mathbb Z$, our induction start. Next, if~\eqref{eq::InductionClaim} holds for a given $k \leq N-1$,~\eqref{eq::y+} yields
\begin{eqnarray}
y_{k+1}^+ - y_{k+2} &\in& (A+BK) \left[  \beta^k (1-\beta) \mathbb Z \right] \notag \\[0.1cm]
&=& \beta^{k+1} (1-\beta) \mathbb Z \; ,
\end{eqnarray}
since $\mathbb Z$ is $\beta$-contractive for the given linear control gain $K$ by construction. This is an induction step implying that~\eqref{eq::InductionClaim} holds for all $k \in \{0,1,\ldots,N\}$. Next, we use the inclusion\footnote{The inclusion $y_{k+1} \in \mathbb Y_{k+1}$ also holds for $k=N-1$, as we have $y_N \in \mathbb T \subseteq Y_N$ due to our particular construction of $\mathbb T$.} $y_{k+1} \in \mathbb Y_{k+1}$ to show that
\begin{eqnarray}
y_{k}^+ &\overset{\eqref{eq::InductionClaim}}{\in}& \mathbb Y_{k+1} \oplus \left[ \beta^k (1-\beta) \mathbb Z \right] \notag \\[0.1cm]
&=& \left[ \mathbb X \ominus (1-\beta^{k+1}) \mathbb Z \right] \oplus \left[ \beta^k (1-\beta) \mathbb Z \right] \notag \\[0.1cm]
&=& \mathbb X \ominus (1-\beta^k) \mathbb Z \; = \; \mathbb Y_{k} ,
\label{eq::YkBound}
\end{eqnarray}
for all $k \in \{ 0,\ldots, N-1 \}$. Similarly, for $k = N$, we have
\begin{eqnarray}
y_{N}^+ &\overset{\eqref{eq::InductionClaim},\eqref{eq::terminalShift}}{\in}& (A+BK)\mathbb T \oplus \left[ \beta^N (1-\beta) \mathbb Z \right] \; \overset{\eqref{eq::Tproperty}}{\subseteq} \; \mathbb T \, . \notag
\end{eqnarray}
Thus, in summary, we have $y_k^+ \in \mathbb Y_k$ for all times indices $k \in \{ 0,1,\ldots, N-1 \}$, $y_N^+ \in \mathbb T$, as well as $w_k^+ \in \mathbb W_k$, where the latter inclusion follows by an argument that is completely analogous to~\eqref{eq::YkBound}. Thus, we have $(y^+,w^+) \in \mathcal A$, which completes the proof.\qed

\section{Real-Time Parallel MPC with Recursive Feasibility Guarantees}
\label{sec::fmpc}

In order to develop a recursively feasible variant of Algorithm~\ref{alg::mpc}, which takes control and state constraints into account, we introduce the auxiliary optimization problem
\begin{align}
V_N(\hat x) \,\defeq \, \underset{x,u}{\min}\;\;& x_N^\tr P x_N+\sum_{k=0}^{N-1} x_k^\tr Q x_k + u_k^\tr R u_k   \notag \\[0.12cm]
\label{eq::MPC2}
\mathrm{s.t.}\;\;&
\left\{
\begin{array}{l}
\forall k \in \{ 0, \ldots, N-1 \}, \\[0.12cm]
x_{k+1} = A x_k + B u_k \; \mid \; \lambda_k , \\[0.12cm]
x_0 = \hat x , \\[0.12cm]
x_{k+1} \in \mathbb Y_{k+1} ,\; u_{k} \in \mathbb W_{k} \; .
\end{array}
\right.
\end{align}
Problems~\eqref{eq::MPC} and~\eqref{eq::MPC2} coincide except for the constraints, which have been replaced in~\eqref{eq::MPC2} by their tightened counterparts, $\mathbb Y_k \subseteq \mathbb X$ and $\mathbb W_k \subseteq \mathbb U$. Consequently,~\eqref{eq::MPC2} can be interpreted as a conservative approximation of~\eqref{eq::MPC}, and we have $V_N \geq J_N$. If $r$ is small, the associated loss of optimality is, however, small, too. We point out that neither~\eqref{eq::MPC} nor~\eqref{eq::MPC2} implement terminal constraints, since $V_N(\hat x) = V_\infty(\hat x)$ as long as $P$ is chosen appropriately and $V_N(\hat x) < \infty$.

\subsection{Recursively feasible parallel MPC}
\label{subsec::recursively_feasible_parallel_MPC}

The main idea of this paper is to solve~\eqref{eq::MPC2} by a variant of Algorithm~\ref{alg::mpc}. For this aim, we introduce the decoupled, and separable, initial problem
\begin{subequations}
\label{eq::initial2}
\begin{align}
\min_{u_0}\;\;&  \Vert u_0 \Vert_R^2 - \lambda_0^\tr B u_0 + \| u_0 - v_0 \|_R^2 
\\[0.12cm]
\mathrm{s.t.} \;\;&
u_0 \in \mathbb U \, , \;
A \hat x + B u_0 \in \mathbb Y_1 \, ,
\end{align}
\end{subequations}
as well as the parametric decoupled QPs
\begin{subequations}
\label{eq::horizon2}
\begin{align}
\underset{x_k,u_k}{\min}\;\; &  \Vert x_k \Vert_Q^2 + \Vert u_k \Vert_R^2
- \lambda_k^\tr B u_k + (\lambda_{k-1} - A^\tr \lambda_k)^\tr x_k  \notag \\[0.12cm]
\qquad \quad&  + \| x_k - z_k \|_Q^2 + \| u_k - v_k \|_R^2
\\[0.12cm]
\mathrm{s.t.}\quad &u_k \in \mathbb W_k \, , \;\;
A x_k + B u_k \in \mathbb Y_{k+1} \, ,
\end{align}
\end{subequations}
defined for all $k \in \{ 1, \ldots, N-1 \}$. The correspond terminal problem~\eqref{eq::terminal} remains unchanged. Finally, a recursively feasible real-time iteration (RFRTI) variant of Algorithm~\ref{alg::mpc} is obtained by replacing Step 1a) with the following Step 1a'):
\begin{itemize}
\item[1a')]  Solve the decoupled QPs~\eqref{eq::initial2},~\eqref{eq::horizon2}, and~\eqref{eq::terminal} and denote their optimal solutions by $x = [x_1,\ldots,x_N]$ and \mbox{$u = [u_0,\ldots,u_{N-1}]$}.
\end{itemize}
From an implementation point of view, this change is minor, as we have merely introduced new (separable) inequality constraints in~\eqref{eq::initial2} and~\eqref{eq::horizon2}, which can still be tackled by explicit MPC solvers~\cite{BorrelliPHD,MPT3}. The advantage of this change is, however, significant as one can now tradeoff conservatism of the constraint margins, controlled by the tuning parameter $r > 0$, with the number $\overline m$ of iterations that are needed to ensure recursive feasibility of the proposed sub-optimal parallel real-time MPC scheme, as the following result holds.

\begin{theorem}
\label{lem::recursive_feasibility}
Let $\sigma > 0$ satisfy $B^T Q B \preceq \sigma R$, let the first state measurement $\hat x$ be such that $V_N(\hat x) < \infty$, and let $P$ and $N$ be chosen such that $V_N(\hat x) = V_\infty(\hat x)$. Then there exists a constants $\kappa < 1$ such that~\eqref{eq::MAIN_INEQUALITY} holds for the iterates of the modified Algorithm 1, where Step~1a) is replaced by Step~1a'). The iterates remain recursively feasible, if $\overline m$ satisfies
\begin{equation}
\label{eq::minimum_mbar}
\overline{m} \geq \log_\kappa{ \left( \frac{1-\beta}{1+\kappa} \cdot \frac{r}{\sigma \Delta(z,v,\lambda)} \right) - 1 } \; ,
\end{equation}
where $\beta < 1$ and $r > 0$ are defined as in Section~\ref{sec::contractive_sets_and_feasibility} and $z,v$, and $\lambda$ denote the initialization of the algorithm.
\end{theorem}

\textit{Proof.} Notice that Theorem~1 in~\cite{Jiang2021} does not use any particular assumption on constraints and, consequently,~\eqref{eq::MAIN_INEQUALITY} also holds for the iterates of the modified Algorithm~\ref{alg::mpc} in the presence of state constraints. Because $V_N(\hat x) = V_\infty(\hat x)$ holds by construction, the terminal region from Lemma~\ref{lem::recursively:admissible} does not have to be implemented and, consequently, Lemma~\ref{lem::recursively:admissible} ensures that recursive feasibility holds whenever the inequality \mbox{$\Vert x_0^+ - x_1^\star \Vert < (1-\beta) r$} is satisfied by the iterates.
Thus, due to Lemma~\ref{lem::distance_to_optimal_finite_iter}, if $\overline m$ satisfies
\[
\Vert x_0^+ - x_1^\star \Vert \leq \sigma (1+\kappa) \kappa^{\overline{m}+1} \Delta(z,v,\lambda) \leq (1 - \beta) r \, .
\]
recursive feasibility is guaranteed. The statement of the theorem follows by solving this inequality with respect to $\overline m$.
\qed

\bigskip
\noindent
Notice that the conditions for asymptotic stability from~\cite{Jiang2021} are unaffected by the state constraints. They are fully applicable to our modified version of Algorithm~1---with the only difference being that we need Theorem~\ref{lem::recursive_feasibility} to guarantee recursive feasibility in the presence of state constraints.

\section{Case Study}
\label{sec::case_study}
Let $p_i$ and $v_i$ denote the position and velocity of the $i$-th cart in a spring-mass-damper chain with $\mathsf{n}$ carts in total. We assume that the system recursion is given by
\begin{eqnarray}
p_i^+ &=& p_i + h v_i \notag \\
v_i^+ &=& v_i + \frac{h}{m} \left( k_\mathrm{s} ( p_{i-1} - 2 p_{i} + p_{i+1} ) -  k_\mathrm{d} v_i + u_i \right) \, \notag
\end{eqnarray}
for all $i \in \{ 1, \ldots, \mathsf{n} \}$, where the mass $m=1$, spring constant $k_\mathrm{s} = 1$, and damping constant $k_\mathrm{d} = 1$ are all set to $1$. Additionally we set $p_{0} = 0$ and $p_{\mathsf{n}+1} = p_\mathsf{n}$ modeling a scenario in which the first cart is attached to a wall while the last cart is free. Here, $u_i$ denotes the force at the $i$-th cart. The discretization parameter is set to $h = 0.1$ while
\begin{eqnarray}
\mathbb X &=& \left\{ x \in \mathbb R^{2 \mathsf{n}} \mid \Vert x \Vert_\infty \leq \overline x \right\} \notag \\ \text{and} \quad \mathbb U &=& \left\{ u \in \mathbb R^{n} \mid \Vert u \Vert_\infty \leq \overline u \right\} \notag
\end{eqnarray}
model symmetric state- and control constraints, where we set $\overline x = \frac{5}{2}$ and $\overline u = 1$. Moreover, we set $Q = I$ and $R = I$; $P$ is found by solving an associated algebraic Riccati equation. In order to study the performance of the proposed algorithm we consider a large-scale MPC problem for $\mathsf{n} = 60$ carts (= $120$ differential states and $60$ control inputs) while the prediction horizon is set to $N = 100$. The initial state measurement is set to $x_0 = 2 \cdot \mathbf{1} \in \mathbb R^{120}$. This leads to a non-trivial sparse QP with $18120$ optimization variables in total, which cannot be solved up to high accuracy in less than $10 \mathrm{s}$ on standard computers.

\begin{figure}
\begin{center}
\includegraphics[scale=0.4]{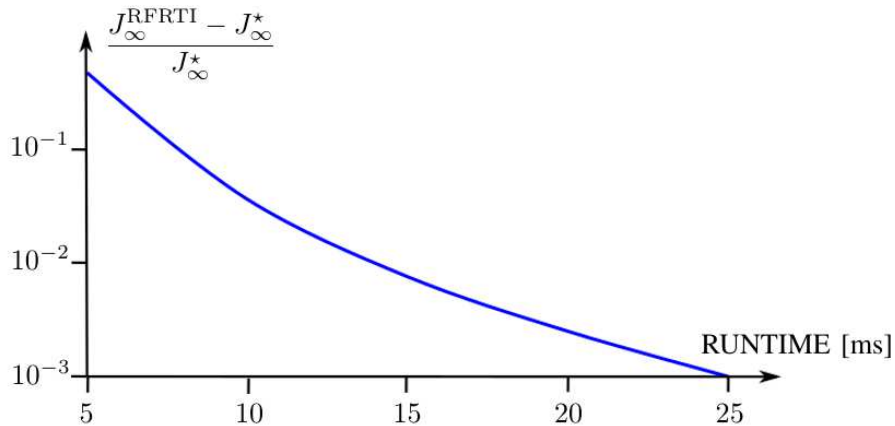}
\end{center}
\caption{\label{fig::loss} The relative loss of performance of RFRTI compared to exact MPC versus the online run-time. Here, $J_\infty^\star$ denotes the (constant) optimal infinite horizon performance obtained by using exact MPC with an online run-time of more than $10 \mathrm{s}$ and $J_\infty^\mathrm{RFRTI}$ denotes the corresponding performance of the proposed real-time controller with recursive feasibility guarantees.}
\end{figure}

\bigskip
\noindent
Next, let $J_\infty^\star$ denote the optimal infinite-horizon closed-loop performance (infinite sum over the stage cost along the closed-loop states and controls) of exact MPC, $J_\infty^\mathrm{RFRTI}$ the measured infinite-horizon closed-loop performance of the proposed RFRTI scheme, and $J_\infty^\mathrm{RTI}$ the corresponding performance of a heuristic real-time iteration without recursive feasibility guarantee. Figure~\ref{fig::loss} shows the relative performance of the proposed RFRTI controller versus its online run-time. Here, we have implemented the proposed algorithm in the form of a prototype Julia code: one inner-loop iteration of the proposed algorithm takes approximately $1 \,$ms; that is, the run-time in milliseconds coincides with $\overline m$. For instance, if we stop each real-time loop after $25 \, \mathrm{ms}$, we obtain a controller that is only $0.1 \%$ suboptimal compared to exact MPC, but a factor $400$ faster, clearly showing the benefit of real-time MPC.

\begin{figure}
\begin{center}
\includegraphics[scale=0.4]{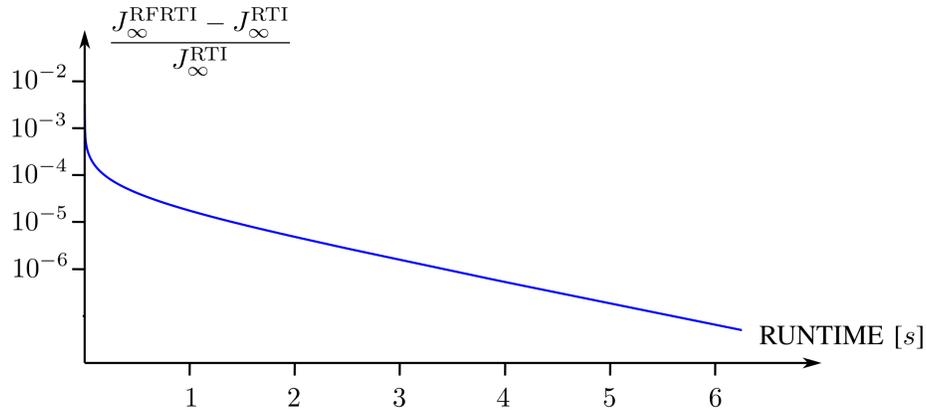}
\end{center}
\caption{\label{fig::RTI} The relative loss of performance, $(J_\infty^\mathrm{RFRTI} - J_\infty^\mathrm{RTI})/J_\infty^\mathrm{RTI}$, where $J_\infty^\mathrm{RFRTI}$ denotes the infinite horizon closed-loop cost of the RTI controller with recursive feasibility guarantee and $J_\infty^\mathrm{RTI}$ the infinite horizon closed-loop cost of a heuristic RTI controller without recursive feasibility guarantee, versus the run-time per-real time iteration.}
\end{figure}

\bigskip
\noindent
Moreover, a heuristic RTI iteration---without the time-varying constraint margins---happens to generate a feasible closed-loop trajectory, although we have no a-priori guarantee for this. Figure~\ref{fig::RTI} compares the relative loss of the RFRTI iteration versus such a heuristic RTI iteration for different online run-times. The relative loss of performance is below $0.1 \%$ in all cases. However, if one wishes to reach a very small loss of performance (e.g. $10^{-6}$), this is only possible if one can accept run-times in the second range. %, which is getting clear from Figure~\ref{fig::RTI}, too.

\section{Conclusions and Outlook}
This paper has presented a RFRTI scheme for model predictive control that comes along with novel recursive feasibility as well as online run-time guarantees (see Lemma~\ref{lem::recursively:admissible} and Theorem~\ref{lem::recursive_feasibility}). A case study for a large-scale linear system with $120$ states and a prediction horizon of $N = 100$ has illustrated the promising performance of the proposed real-time controller compared to exact MPC. Moreover, it is has been found that enforcing recursive feasibility guarantees comes at a negligible loss of performance. Future research will focus on an open-source software implementation of the proposed RFRTI control scheme.


\begin{thebibliography}{00}

\bibitem{Rawlings2009}
J.~Rawlings, D.~Mayne, and M.~Diehl, {\em Model Predictive Control: Theory and
  Design}.
\newblock Madison, WI: Nob Hill Publishing, 2018.

\bibitem{Qin2003}
S.~Qin and T.~Badgwell, ``A survey of industrial model predictive control
  technology,'' {\em Control Engineering Practice}, vol.~93, no.~316,
  pp.~733--764, 2003.

\bibitem{Ferreau2014}
H.~J. Ferreau, C.~Kirches, A.~Potschka, H.~G. Bock, and M.~Diehl, ``qpoases: A
  parametric active-set algorithm for quadratic programming,'' {\em
  Mathematical Programming Computation}, vol.~6, no.~4, pp.~327--363, 2014.

\bibitem{FRISON2020}
G.~Frison and M.~Diehl, ``Hpipm: a high-performance quadratic programming
  framework for model predictive control,'' {\em IFAC-PapersOnLine}, vol.~53,
  no.~2, pp.~6563--6569, 2020.
\newblock 21st IFAC World Congress.

\bibitem{Everett1963}
H.~Everett, ``Generalized {L}agrange multiplier method for solving problems of
  optimum allocation of resources,'' {\em Operations Research}, vol.~11, no.~3,
  pp.~399--417, 1963.

\bibitem{Boyd2011}
S.~Boyd, N.~Parikh, E.~Chu, B.~Peleato, and J.~Eckstein, ``Distributed
  optimization and statistical learning via the alternating direction method of
  multipliers,'' {\em Foundation Trends in Machine Learning}, vol.~3, no.~1,
  pp.~1--122, 2011.

\bibitem{Houska2016}
B.~Houska, J.~Frasch, and M.~Diehl, ``An augmented {L}agrangian based algorithm
  for distributed non-convex optimization,'' {\em SIAM Journal on
  Optimization}, vol.~26, no.~2, pp.~1101--1127, 2016.

\bibitem{Giselesson2013}
P.~Giselsson, M.~Dang~Doan, T.~Keviczky, B.~De~Schutter, and A.~Rantzer,
  ``Accelerated gradient methods and dual decomposition in distributed model
  predictive control,'' {\em Automatica}, vol.~49, no.~3, pp.~829--833, 2013.

\bibitem{Conte2012}
C.~Conte, T.~Summers, M.~Zeilinger, M.~Morari, and C.~Jones, ``Computational
  aspects of distributed optimization in model predictive control,'' in {\em
  Proceedings of the 51st IEEE Conference on Decision and Control, 2012},
  pp.~6819--6824, 2012.

\bibitem{Necoara2008}
I.~Necoara and J.~Suykens, ``Application of a smoothing technique to
  decomposition in convex optimization,'' {\em IEEE Transactions on Automatic
  Control}, vol.~53, no.~11, pp.~2674--2679, 2008.

\bibitem{Donoghue2013}
B.~O'Donoghue, G.~Stathopoulos, and S.~Boyd, ``A splitting method for optimal
  control,'' {\em IEEE Transactions on Control Systems Technology}, vol.~21(6),
  pp.~2432--2442, 2013.

\bibitem{Jiang2021}
Y.~Jiang, J.~Oravec, B.~Houska, and M.~Kvasnica, ``Parallel {MPC} for linear
  systems with input constraints,'' {\em IEEE Transactions on Automatic
  Control}, vol.~66, no.~7, pp.~3401--3408, 2021.

\bibitem{BemEtal:aut:02}
A.~Bemporad, M.~Morari, V.~Dua, and E.~Pistikopoulos, ``The explicit linear
  quadratic regulator for constrained systems,'' {\em Automatica}, vol.~38,
  no.~1, pp.~3--20, 2002.

\bibitem{Ingole2015}
D.~Ingole and M.~Kvasnica, ``{FPGA} implementation of explicit model predictive
  control for closed loop control of depth of anesthesia,'' in {\em 5th IFAC
  Conference on Nonlinear Model Predictive Control}, pp.~484--489, 2015.

\bibitem{borrelli2003geometric}
F.~Borrelli, A.~Bemporad, and M.~Morari, ``Geometric algorithm for
  multiparametric linear programming,'' {\em Journal of optimization theory and
  applications}, vol.~118, no.~3, pp.~515--540, 2003.

\bibitem{FERRAMOSCA2013}
A.~Ferramosca, D.~Limon, I.~Alvarado, and E.~Camacho, ``Cooperative distributed
  {MPC} for tracking,'' {\em Automatica}, vol.~49, no.~4, pp.~906--914, 2013.

\bibitem{Conte2016}
C.~Conte, C.~Jones, M.~Morari, and M.~Zeilinger, ``Distributed synthesis and
  stability of cooperative distributed model predictive control for linear
  systems,'' {\em Automatica}, vol.~69, pp.~117--125, 2016.

\bibitem{Darivianakis::TAC19}
G.~Darivianakis, A.~Eichler, and J.~Lygeros, ``Distributed model predictive
  control for linear systems with adaptive terminal sets,'' {\em IEEE
  Transactions on Automatic Control}, vol.~65, no.~3, pp.~1044--1056, 2020.

\bibitem{HERNANDEZ201711829}
B.~Hernandez, P.~Baldivieso, and P.~Trodden, ``Distributed {MPC}: Guaranteeing
  global stability from locally designed tubes,'' {\em IFAC-PapersOnLine},
  vol.~50, no.~1, pp.~11829--11834, 2017.
\newblock 20th IFAC World Congress.

\bibitem{Cannon::TAC::2017::7488973}
M.~Schulze~Darup and M.~Cannon, ``On the computation of $\lambda$- contractive
  sets for linear constrained systems,'' {\em IEEE Transactions on Automatic
  Control}, vol.~62, no.~3, pp.~1498--1504, 2017.

\bibitem{Matthias17}
M.~A. M\"uller and F.~Allg\"ower, ``Economic and distributed model predictive
  control: Recent developments in optimization-based control,'' {\em SICE
  Journal of Control, Measurement, and System Integration}, vol.~10, no.~2,
  pp.~39--52, 2017.

\bibitem{Damm2014}
T.~Damm, L.~Gr\"une, M.~Stieler, and K.~Worthmann, ``An Exponential Turnpike Theorem for Dissipative Discrete Time Optimal Control Problems'' {\em SIAM Journal on Control and Optimization}, vol.~52, no.~3, pp.~1935--1957, 2014.

\bibitem{BorrelliPHD}
F.~Borrelli, {\em Constrained Optimal Control Of Linear And Hybrid Systems},
  vol.~290 of {\em Lecture Notes in Control and Information Sciences}.
\newblock Springer, 2003.

\bibitem{MPT3}
M.~Herceg, M.~Kvasnica, C.~Jones, and M.~Morari, ``Multi-parametric toolbox
  3.0,'' in {\em 2013 European Control Conference}, pp.~502--510, 2013.

\bibitem{Bertsekas2012}
D.~Bertsekas, {\em Dynamic Programming and Optimal Control}.
\newblock Belmont, Massachusetts: Athena Scientific Dynamic Programming and
  Optimal Control, 3rd~ed., 2012.

\bibitem{Blanchini2008}
F.~Blanchini and S.~Miani, {\em Set-theoretic methods in control}.
\newblock Systems \& Control: Foundations \& Applications, Birkh\"auser Boston,
  Inc., Boston, MA, 2008.

\end{thebibliography}
\end{document}